\documentclass[1p]{elsarticle}
\usepackage{amssymb}
\usepackage{amsfonts}

\newtheorem{theorem}{Theorem}
\newtheorem{conjecture}[theorem]{Conjecture}

\newproof{pf}{Proof}

\begin{document}

\title{The 1--2--3 Conjecture holds for graphs with large enough minimum degree}

\author[agh]{Jakub Przyby{\l}o} 


\address[agh]{AGH University of Science and Technology, al. A. Mickiewicza 30, 30-059 Krakow, Poland}

\begin{abstract}
A simple graph more often than not contains adjacent vertices with equal degrees. This in particular holds for all pairs of neighbours in regular graphs, while a lot such pairs can be expected  e.g. in many random models. Is there a universal constant $K$, say $K=3$, such that one may always dispose of such pairs from any given connected graph with at least three vertices by blowing its selected edges into at most $K$ parallel edges? This question was first posed in 2004 by Karo\'nski,  \L uczak and Thomason, who equivalently asked if one may assign weights $1,2,3$ to the edges of every such graph so that adjacent vertices receive distinct weighted degrees -- the sums of their incident weights. This basic problem is commonly referred to as the 1--2--3 Conjecture nowadays, and has been  addressed in multiple papers. Thus far it is known that weights $1,2,3,4,5$ are sufficient [J. Combin. Theory Ser. B 100 (2010) 347--349]. We show that this conjecture holds if only the minimum degree $\delta$ of a graph is large enough, i.e. when $\delta = \Omega(\log\Delta)$, where $\Delta$ denotes the maximum degree of the graph. The principle idea behind our probabilistic proof relies on associating random variables with a special and carefully designed distribution to most of the vertices of a given graph, and then choosing weights for major part of the edges depending on the values of these variables in a deterministic or random manner.  
\end{abstract}


\maketitle

\section{introduction}
Let $G=(V,E)$ be a (simple) graph. A \emph{$k$-weighting} of $G$ is a mapping $\omega:E\to [k]$, where $[k]=\{1,2,\ldots,k\}$. Given such $\omega$ and $v\in V$, we denote by $E(v)$ the set of edges incident with $v$ in $G$ and by $s_\omega(v):=\sum_{e\in E(v)}\omega(e)$ the so-called \emph{weighted degree} of the vertex $v$, abbreviating its notation to $s(v)$ if this causes no ambiguities, and refer to it  simply as the \emph{sum of $v$}. We say that vertices $u,v$ are \emph{in sum conflict} if $uv\in E$ and $s(u)=s(v)$. If there are no sum conflicts, i.e. when $s_\omega$ is a proper vertex colouring of $G$, we moreover say that $\omega$ is \emph{vertex-colouring}. In 2004 Karo\'nski, {\L}uczak and Thomason~\cite{123KLT} asked a basic question referred to ever since as the \emph{1--2--3 Conjecture}:
\begin{conjecture}\label{Conjecture123Conjecture}
Every connected graph with at least two edges admits a ver\-tex-co\-louring $3$-weighting.
\end{conjecture}

One of the main motivations behind this problem are research devoted to irregularities in graphs. 
A given graph or multigraph is called \emph{irregular} (cf.~\cite{Chartrand}) if all its vertices have pairwise distinct degrees, whereas we call it \emph{locally irregular} (cf.~\cite{LocalIrreg_1,BensmailMerkerThomassen,LocalIrreg_2}) if the degree of its every vertex is distinct from the degrees of all its neighbours.  
While regular graphs -- graphs with equal degrees of all vertices are fundamental and thoroughly investigated objects in graph theory, their antonyms do not even exist, as a basic application of the pigeonhole principle implies there are no irregular graphs with more than one vertex. In~\cite{ChartrandErdosOellermann} Chartrand, Erd\H{o}s and Oellermann thus proposed and studied several alternative definitions of irregular graphs, while in~\cite{Chartrand} Chartrand et al. 
initiated a related study of a graph invariant designed to capture in some sense the level of irregularity of a graph.
Namely, they defined the \emph{irregularity strength} of a graph $G=(V,E)$ as the least $k$ such that (similarly as above) a $k$-weighting of $G$ exists inducing distinct weighted degrees of \emph{all} vertices (regardless of their adjacency), i.e. such that $s(u)\neq s(v)$ whenever $u,v\in V$ and $u\neq v$.  This concept exploits the fact that irregular multigraphs, unlike irregular graphs are pretty common objects, for the irregularity strength of a (simple) graph $G$ can equivalently be defined as the least $k$ for which we may create an irregular multigraph of $G$ by blowing  each edge $e$ of $G$ up to at most $k$ parallel copies of $e$ 
(where multiplicities of edges correspond to their weights in the first definition). 
In the same vein, Conjecture~\ref{Conjecture123Conjecture} asserts that we may obtain a \emph{locally} irregular multigraph of every connected graph $G$ of size at least two by blowing each edge $e$ of $G$ up to at most three copies of $e$.

Numerous papers have been devoted to irregularity strength of graphs. 
Many of these, e.g.~\cite{Amar,Lazebnik,Faudree,Frieze,KalKarPf,MajerskiPrzybylo2,AsymptoticIrregStrReg,Przybylo}
concentrate around perhaps the most well-known problem related with this graph invariant, concerning a general upper bound within the class of regular graphs and a conjecture of Faudree and Lehel~\cite{Faudree} 
(first posed in the form of question by Jacobson, as mentioned by Lehel in~\cite{Lehel}),
which remains open since~1987. See also e.g.~\cite{Aigner,Bohman_Kravitz,Dinitz,Ebert,Faudree2,Ferrara,Frieze,Gyarfas,KalKarPf,MajerskiPrzybylo2,Nierhoff,irreg_str2} for other exemplary results.
A large number of related concepts and intriguing open questions were inspired by this problem.
The irregularity strength was in fact the cornerstone of  nowadays fast-developing wide  branch of research on various variants of graph or hypergraph colourings and labellings, see~\cite{Gallian_survey,Lehel,Seamon123survey} for surveys devoted to some of them.
One of these problems, concerning proper edge colourings associating distinct sets of incident colours to adjacent vertices, see e.g.~\cite{BalGLS,Hatami,Joret,Zhang} for details (or~\cite{Louigi2,123KLT,Vuckovic_3-multisets} for its non-proper counterpart), was mentioned already in the first paper~\cite{123KLT} on the 1--2--3 Conjecture, and can be regarded as one of its main inspirations, apart from its obvious progenitor -- the irregularity strength. 

The concept of vertex-colouring weightings gained equally considerable attention in the combinatorial community as its precursor largely due to the mentioned beautiful Conjecture~\ref{Conjecture123Conjecture}, which has a very simple and elegant formulation and no obvious solution. The introductory paper~\cite{123KLT} devoted to it provided no constant upper bound in the general case (this was only assured for the variant of the problem admitting real numbers as weights), but settled validity of Conjecture~\ref{Conjecture123Conjecture} for the family of graphs with the chromatic number $\chi(G)$ at most $3$. The first finite upper bound was exhibited by Addario-Berry, Dalal, McDiarmid, Reed and Thomason~\cite{Louigi30}, who proved that $K=30$ assures existence of  a ver\-tex-co\-louring $K$-weighting for every connected graph with at least two edges. This constant was pushed down to $K=16$ by Addario-Berry,  Dalal and Reed~\cite{Louigi}, and then to $K=13$ by Wang and Yu~\cite{123with13} by means of a similar approach as in~~\cite{Louigi30}, based on special theorems on so-called degree constrained subgraphs. In~\cite{Louigi} it was moreover observed that even just weights $1,2$ are asymptotically almost surely sufficient for a random graph (chosen from $G_{n,p}$ for a constant $p\in(0,1)$). A big breakthrough was later achieved due to development of a simple algorithm by Kalkowski~\cite{Kalkowski12} (designed to tackle a total variant of the 1--2--3 Conjecture, see~\cite{12Conjecture}), whose refinement allowed Kalkowski, Karo\'nski and Pfender~\cite{KalKarPf_123} to narrow down the set of necessary weights to merely $\{1,2,3,4,5\}$. Recently Przyby{\l}o~\cite{1234Reg123} proved that the weight $5$ is redundant in the case of regular graphs (see~\cite{Julien5regular123} for an earlier results in the case of $5$-regular graphs)
and confirmed the conjecture for $d$-regular graphs with $d$ large enough. Lately also a certain family of very dense graphs for which weights $1,2,3$ suffice was exposed by Zhong~\cite{123dense-Zhong}, who showed  there exists a constant $n'$ such that every graph with $n\geq n'$ vertices and minimum degree $\delta(G)>0.99985n$ admits a vertex-colouring edge $3$-weighting. On the other hand it is known that graphs which require weights $1,2,3$ are not so uncommon, as confirmed by Dudek and Wajc~\cite{DudekWajc123complexity}, who proved that determining whether a particular graph admits a vertex-colouring $2$-weighting is NP-complete. This however ceases to hold in the bipartite case due to relatively simple description of those bipartite graphs which necessitate use of weight $3$ provided by  Thomassen, Wu and Zhang~\cite{ThoWuZha}, which implies a polynomial time algorithm deciding when this is the case. Similarly as the irregularity strength of graphs, the 1--2--3 Conjecture also gave rise to a list of intriguing related concepts, as e.g. its extension to a natural list setting~\cite{BarGrNiw,PrzybyloWozniakChoos,WongZhuChoos}, with interesting applications~\cite{BarGrNiw,WongZhu23Choos} of algebraic approach exploiting Alon's Combinatorial Nullstellensatz. Intriguingly, though it is believed that 3-element lists of weights should suffice, see~\cite{BarGrNiw}, no finite upper bound is known thus far in this more demanding setting; see~\cite{Logarithmic_Weight-choosability} for a result implying that lists of length $O(\log\Delta)$ are enough. 

As the main contribution of this paper we confirm that 1--2--3 Conjecture indeed holds if only the minimum degree $\delta$ of a graph is larger than a constant times the logarithm of its maximum degree $\Delta$ (i.e. $\delta=\Omega(\log\Delta)$). 

\begin{theorem}\label{123ConjectureForLargedelta}
There exists a constant $C$ such that every graph with maximum degree $\Delta\geq 2$ and minimum degree $\delta\geq C \log\Delta$  admits a ver\-tex-co\-louring $3$-weighting.
\end{theorem}
Here, by `$\log$' we mean the natural logarithm.
We shall be further on using `$\ln$' instead in order to shorten notation.

\section{Outline of Main Ideas Exploited in the Proof of Theorem~\ref{123ConjectureForLargedelta}}

Suppose we are given a graph $G=(V,E)$, and we wish to use the probabilistic method to design its desired weighting using mostly weights $1$ and $3$. 

The very general yet still vague main idea of the proof relies on assigning random variables $X_v$ to the vertices $v\in V$ and then associating weight $3$ to any given edge $uv$ if $X_u+X_v$ is large enough, and assigning weight $1$ to it otherwise. There are several deviations from this general rule though.

Suppose first that we associate with every vertex $v$ a random variable with uniform distribution over the interval $[1,3]$. We then assign weight $3$ to $uv\in E$ if $X_u+X_v\geq 4$, and weight $1$ otherwise. Then, for each $x\in[1,3]$, we have:
$\mathbf{E}(s(v)|X_v=x)=\sum_{u\in N(v)}(3\cdot \mathbf{Pr}(X_u\geq 4-x)+1\cdot \mathbf{Pr}(X_u< 4-x)) = d(v)(3\cdot \frac{x-1}{2} + 1\cdot \frac{3-x}{2})=xd(v)$. Such general idea yields a convenient initial distribution of sums e.g. for $d$-regular graphs, as every vertex may then receive sums in the range $[d,3d]$, which tend to be fairly evenly distributed throughout the graph, and at most $d$ of these  might be `blocked' by its neighbours. However, such approach fails in case of non-regular graphs, for suppose all $d(v)$ neighbours of $v$ have e.g. degree $0.4d(v)$ and $X_v\in[1,1.1]$.
Then the expected number of neighbours $u$ of $v$ with sums in $[d(v),1.1d(v)]$ shall roughly equal $d(v)/8>|[d(v),1.1d(v)]|$.

In order to counteract against such a potential problem, we shall thus use a two-fold modified random approach. First of all, random variables $X_v$ assigned to the vertices of $G$ shall have distribution with constant density replaced by a function $g(x)=\frac{c}{x}$ for appropriately chosen $c$. Thereby, we shall prevent high concentration of sums of the neighbours of any vertex $v$ with neighbourhood consisting of many vertices of `slightly' smaller degrees than $d(v)$ (the price to pay for this is a larger probability of a vertex to have smaller sum, i.e. closer to $d(v)$ than $3d(v)$ -- the choice of density function is optimized with respect to this phenomenon though). The second major and necessary alteration of the random process concerns the rules of assigning weights to the edges within it. In many cases we shall still assign weight $3$ to a given edge $uv$ if $X_u+X_v$ is large enough, but this time `large enough' shall be expressed by a kind of two-variable function, designed carefully so that we actually achieve as a result of the process the assumed desired distribution of sums of the vertices (with roughly $\mathbf{E}(s(v)|X_v=x)=xd(v)$, where the sum of $v$ is in particular more likely to be smaller than larger). Moreover, in order to achieve exactly such a goal, for some values of $X_u$ and $X_v$ we shall choose weight $1$ or $3$ for  $uv\in E$ randomly, with a given probability dependent on the values of $X_u$ and $X_v$ -- this shall be in practice realized via prior assignment of independent random variables $X_e$ with uniform distribution over $[0,1]$ to (almost) all edges $e\in E$, and then using their values for comparison with a predefined function of $X_u$ and $X_v$, if needed, in order to chose the weight for $uv$. 

Furthermore, so that all requirements above could be fulfilled, we must use as domains of our random variables $X_v$ the set $[1.1,2.9]$, not $[1,3]$ (for which this is not possible), thus the sum of (almost) all vertices $v$ shall belong to: $[(1.1-\varepsilon)d(v),(2.9+\varepsilon)d(v)]$.

As we shall heavily rely on probability, we shall be able to control only the approximate size of most of the sums, thus, in order to finally adjust them to some precise values, we shall beforehand (also randomly) ration out from $V$ some relatively small subset $U$ and perform our random edge weight assignment within the set $W:=V\smallsetminus U$ (containing vast majority of all the vertices and edges of $G$). Weights of a part of the edges between $U$ and $W$, denoted by $F$, shall then be used to adjust the sums in $U$ and $W$. We shall thus set apart their (relatively small) subset $F_W$ to be later utilized to handle the sums in $W$, and $F_U$ -- useful to distinguish the vertices in $U$. The choice of weights for edges in $F_U$ shall guarantee that only a limited number of neighbours in $U$ of any vertex in $U$ might ultimately have the same sum as this vertex, small enough that a special algorithm shall eventually allow us to distinguish all vertices in $U$, even though we shall admit assigning them only sums congruent to $0$ or $1$ modulo $100$. These residues modulo $100$ shall be forbidden on the other hand for the sums in $W$, thus admitting no sum conflicts between vertices in $U$ and $W$; all details follow.

\section{Tools}

For random arguments we shall in particular need the symmetric variant of the Lov\'asz Local Lemma, see e.g.~\cite{AlonSpencer} and the Chernoff Bound, see e.g.~\cite{JansonLuczakRucinski}. 

\begin{theorem}[\textbf{The Local Lemma}]
\label{LLL-symmetric}
Let $A_1,A_2,\ldots,A_n$ be events in an arbitrary pro\-ba\-bi\-li\-ty space.
Suppose that each event $A_i$ is mutually independent of a set of all the other
events $A_j$ but at most $D$, and that $\mathbf{Pr}(A_i)\leq p$ for all $1\leq i \leq n$. If
$$ p \leq \frac{1}{e(D+1)},$$
then $ \mathbf{Pr}\left(\bigcap_{i=1}^n\overline{A_i}\right)>0$.
\end{theorem}
\begin{theorem}[\textbf{Chernoff Bound}]\label{ChernofBoundTh}
For any $0\leq t\leq np$,
$$\mathbf{Pr}\left(\left|{\rm BIN}(n,p)-np\right|>t\right)<2e^{-\frac{t^2}{3np}}$$
where ${\rm BIN}(n,p)$ is the sum of $n$ independent Bernoulli variables, each equal to $1$ with probability $p$ and $0$ otherwise.
\end{theorem}

\section{Proof of Theorem~\ref{123ConjectureForLargedelta}}

It is sufficient to prove the theorem for graphs with maximum degree $\Delta$ large enough.
Let $G=(V,E)$ be a graph with $\delta=\delta(G)\geq 10^{20}\ln\Delta$ and  $\Delta\geq 2$ sufficiently large so that all explicit inequalities below hold (in particular, $G$ has no isolated edges).

We shall first ration out in Subsection~\ref{Subsection_U_FW} an auxiliary set $U$ from $V$ and a special subset $F_W$ of the edges between $U$ and $W=V\smallsetminus U$, whose weights shall be later used to adjust the sums in $W$. In the subsequent Subsection~\ref{Subsection_FU} we shall in turn single out a special subset $F_U$ of the remaining edges between $U$ and $W$ with several features, which shall in particular help us scatter sums of the majority of neighbouring vertices in $U$ (whose most incident edges shall join them with $W$) relatively far apart from each other -- these shall be ultimately sum distinguished within the algorithm in the last Subsection~\ref{Subsection_Ufinal}. Subsection~\ref{Subsection_Wpreliminaries} in turn includes preliminary definitions and quantities essential for the main random procedure (exploiting random variables $X_v$ and $X_e$) of weighting the edges incident with vertices in $W$, carried out in Subsections~\ref{Subsection_Rules_W} -- \ref{Subsection_Wfinal}.

By $d(v)$ and $N(v)$ we shall understand below the degree and the neighbourhood of $v\in V$ in $G$, i.e. $d_G(v)$ and $N_G(v)$, respectively. Moreover, for subsets of vertices $A,B\subseteq V$, by $d_A(v)$ we mean the number of neighbours of $v$ (in $G$) belonging to $A$, while by $E(A,B)$ we denote the set of edges $e\in E$ with one end in $A$ and the other one in $B$, and by $E(A)$ -- the set of edges $e\in E$ with both ends in $A$.

\subsection{Choosing $U$ and $F_W$}\label{Subsection_U_FW}

We first construct $U\subseteq V$ such that for every $v\in V$:
\begin{equation}\label{Udegree}
|d_U(v)-10^{-4}d(v)|\leq 10^{-6}d(v).
\end{equation}
Suppose we choose every vertex $v\in V$ to be placed in $U$ randomly and independently with probability $10^{-4}$
and let us consider the following event:
$$A_{1,v}~:~ |d_U(v)-10^{-4}d(v)|>10^{-6}d(v).$$
By the Chernoff Bound (Theorem~\ref{ChernofBoundTh}),
$$\mathbf{Pr}\left(A_{1,v}\right)\leq 2e^{-\frac{\left(10^{-6}d(v)\right)^2}{3\cdot 10^{-4}d(v)}}
= 2e^{-\frac{d(v)}{3\cdot 10^{8}}}
\leq 2e^{-\frac{10^{20}\ln\Delta}{3\cdot 10^{8}}}
= \frac{2}{\Delta^{\frac{10^{12}}{3}}}
\leq \frac{1}{e\Delta^2}.$$
As every event $A_{1,v}$ is mutually independent of all other events $A_{1,u}$ except these with $N(u)\cap N(v)\neq \emptyset$, i.e. all except at most $\Delta (\Delta-1)$, by the Lov\'asz Local Lemma (Theorem~\ref{LLL-symmetric}), with positive probability none of the events $A_{1,v}$ appears, and thus we may choose $U\subseteq V$ so that (\ref{Udegree}) is fulfilled for every $v\in V$. 

We fix any such $U$ and set 
\begin{equation}\label{WandFDefinition}
W:=V\smallsetminus U~~~~~~~~{\rm and}~~~~~~~~ F:=E(U,W).
\end{equation}

We shall next choose $F_W\subseteq F$ so that for every $w\in W$ and each $u\in U$:
\begin{equation}\label{F_WdegreeInW}
|d_{F_W}(w)-10^{-4}d_U(w)|\leq 10^{-6}d_U(w),
\end{equation}
\begin{equation}\label{F_WdegreeInU}
|d_{F_W}(u)-10^{-4}d_W(u)|\leq 10^{-6}d_W(u).
\end{equation}
Analogously as above, we place every edge $e\in F$ in $F_W$ randomly and independently with probability $10^{-4}$. 
We then denote the following events for $w\in W$ and $u\in U$:
$$
A_{2,w}~:~|d_{F_W}(w)-10^{-4}d_U(w)| > 10^{-6}d_U(w),
$$
$$
A_{3,u}~:~|d_{F_W}(u)-10^{-4}d_W(u)|> 10^{-6}d_W(u).
$$
By the Chernoff Bound and (\ref{Udegree}):
\begin{equation}\label{PrA2w}
\mathbf{Pr}\left(A_{2,w}\right) < 2e^{-\frac{d_U(w)}{3\cdot 10^8}} 
\leq 2e^{-\frac{\left(10^{-4}-10^{-6}\right)d(w)}{3\cdot 10^8}} 
\leq 2e^{-\frac{10^{20}\ln\Delta}{4\cdot 10^{12}}} 
\leq \frac{1}{e\Delta},
\end{equation}
\begin{equation}\label{PrA3w}
\mathbf{Pr}\left(A_{3,u}\right) < 2e^{-\frac{d_W(u)}{3\cdot 10^8}} 
\leq 2e^{-\frac{\left(1-10^{-4}-10^{-6}\right)d(u)}{3\cdot 10^8}} 
\leq 2e^{-\frac{10^{20}\ln\Delta}{4\cdot 10^8}} 
\leq \frac{1}{e\Delta}.
\end{equation}
Every event $A_{2,w}$ is mutually independent of all other events $A_{2,v}$, $A_{3,v}$ except at most $d_U(w)<d(w)$ of them, while every $A_{3,u}$ is mutually independent of all other events $A_{2,v}$, $A_{3,v}$ except at most $d_W(u)<d(w)$. Therefore, the Lov\'asz Local Lemma, (\ref{PrA2w}) and~(\ref{PrA3w}) imply that with positive probability none of the events  $A_{2,v}$, $A_{3,v}$ appears, and hence there is $F_W\subseteq F$ consistent with (\ref{F_WdegreeInW}) and (\ref{F_WdegreeInU}).

Fix any such $F_W$ and set: 
\begin{equation}\label{F'_definition}
F':=F\smallsetminus F_W.
\end{equation}

\subsection{Designing $F_U$}\label{Subsection_FU}

In order to obtain a convenient distribution of sums within $U$ in a further part of our construction, not influencing at the same time differentiation of the sums in $W$ significantly, we design $F_U\subseteq F'$ in two steps (resulting in diversified expectations of $d_{F_U}(u)$ for vertices $u\in U$, and unified in terms of $d_{F'}(w)$ for $w\in W$, cf.~(\ref{FUdegreeInU}) and~(\ref{FUdegreeInW})).

First for every vertex  $u\in U$ we independently and equiprobably choose an integer $i_u\in \{0,1,\ldots,10^3-1\}$, and then we independently place every edge $uw\in F'$ with $u\in U$ in $F_U$ with probability $10^{-3}i_u$. We shall prove that $F_U\subseteq F'$ (and $i_u$, $u\in U$) can be chosen so that, in particular, for every $u\in U$ and $w\in W$:  
\begin{eqnarray}
|d_{F_U}(u) - 10^{-3}i_ud_{F'}(u)| &\leq& 10^{-5}d(u), \label{FUdegreeInU}\\
|d_{F_U}(w) - \frac{1-10^{-3}}{2}d_{F'}(w)|&\leq& 10^{-5}d_{F'}(w). \label{FUdegreeInW}
\end{eqnarray}
For this aim, let us denote the following events for $u\in U$ and $w\in W$:  
\begin{eqnarray}
&&A_{4,u}~:~ |d_{F_U}(u) - 10^{-3}i_ud_{F'}(u)| > 10^{-5}d(u), \nonumber\\
&&A_{5,w}~:~ |d_{F_U}(w) - \frac{1-10^{-3}}{2}d_{F'}(w)| > 10^{-5}d_{F'}(w). \nonumber
\end{eqnarray} 
By the Chernoff Bound, for every $u\in U$ and every fixed $i\in\{1,2,\ldots,10^3-1\}$:
$$
\mathbf{Pr}\left(A_{4,u}~|~i_u=i\right) < 2e^{-\frac{\left(10^{-5}d(u)\right)^2}{3\cdot 10^{-3}i\cdot d_{F'}(u)}} 
\leq 2e^{-\frac{\left(10^{-5}d(u)\right)^2}{3\cdot d(u)}} 
\leq 2e^{-\frac{10^{10}\ln\Delta}{3}} <\frac{1}{e\cdot 3\Delta^2},
$$
while we obviously have: $\mathbf{Pr}\left(A_{4,u}~|~i_u=0\right) = 0$. 
Thus, by the low of total probability:
\begin{equation}\label{PrA4u}
\mathbf{Pr}\left(A_{4,u}\right) <\frac{1}{e\cdot 3\Delta^2}.
\end{equation}
Note further that for any given vertex $w\in W$ and an edge $uw\in F'$, the probability that $uw$ is placed in $F_U$ equals $10^{-3}\cdot[0\cdot10^{-3}+1\cdot10^{-3}+\ldots+(10^3-1)\cdot10^{-3}] = \frac{1-10^{-3}}{2}$ and such events are independent for distinct edges incident with a fixed $w\in W$. Thus, by the Chernoff Bound, (\ref{Udegree}) and (\ref{F_WdegreeInW}): 
\begin{equation}\label{PrA5w}
\mathbf{Pr}\left(A_{5,w}\right) < 2e^{-\frac{\left(10^{-5}d_{F'}(w)\right)^2}{3\cdot \frac{1-10^{-3}}{2}d_{F'}(w)}}
\leq  2e^{-\frac{10^{-10}d_{F'}(w)}{1.5}}
\leq 2 e^{-\frac{10^{-10}\cdot 0.5 \cdot 10^{-4}d(w)}{1.5}}
<\frac{1}{e\cdot 3\Delta^2}.
\end{equation}
 
 We shall require one more feature to be obeyed by the choices of $F_U$ and $i_u$, $u\in U$ (accountable for a later convenient sum distribution in $U$). In order to phrase it, let us first  define for every $u\in U$ an interval:
 \begin{eqnarray}
 J(u):&=& [d(u)+10^{-3}i_ud_{F'}(u)-10^{-5}d(u),\nonumber\\
        & & ~d(u)+10^{-3}i_ud_{F'}(u)+10^{-5}d(u)+d_{F_W}(u) + 2d_U(u)] \label{JuDefinition}
 \end{eqnarray}
(which eventually shall contain the final sum of $u$). 
Note that by (\ref{F_WdegreeInU}) and (\ref{Udegree}), one may bound the length of $J(u)$ as follows: 
 \begin{eqnarray}
 |J(u)| &=& 2\cdot10^{-5}d(u)+d_{F_W}(u) + 2d_U(u)\nonumber\\ 
 &\leq& 2\cdot10^{-5}d(u) + (10^{-4}+10^{-6})d_W(u) + 2(10^{-4}+10^{-6})d(u) \nonumber\\
&<& 3.23\cdot 10^{-4} d(u), \label{JuSize}
 \end{eqnarray}
 while the distance between the left ends of two consecutive such possible intervals (i.e. obtained for $i_u=i$ and $i_u=i+1$ for some integer $i$), by~(\ref{F'_definition}), (\ref{F_WdegreeInU}) and (\ref{Udegree}), equals at least:
 \begin{eqnarray}
 10^{-3}d_{F'}(u) &\geq&
 10^{-3}(d_W(u)-d_{F_W}(u)) \geq 10^{-3}(1-10^{-4}-10^{-6})d_W(u) \nonumber\\
 &\geq& 10^{-3}(1-10^{-4}-10^{-6})^2 d(u)
> 9.9 \cdot 10^{-4}d(u). \label{3TimesJuLength}
\end{eqnarray}
As at the end of our construction we shall need to distinguish any given vertex $u\in U$ only from its neighbours $u'\in U$ with $0.5d(u)\leq d(u')\leq d(u)$, while the sum of every vertex $v\in U$ shall belong to $J(v)$, the only vertices which might end up in sum conflict with a given $u\in U$ shall belong to the following set:
\begin{equation}\label{NULeqDef}
N^{U}_\leq(u):=\{u'\in N_U(u):~0.5d(u)\leq d(u')\leq d(u)~\wedge~J(u')\cap J(u)\neq \emptyset\}.
\end{equation}
We shall thus require that apart from~(\ref{FUdegreeInU}) and~(\ref{FUdegreeInW}) also the following shall hold for every $u\in U$:
\begin{equation}\label{NULeqInU}
|N^{U}_\leq(u)| \leq 2\cdot 10^{-3}d_U(u).
\end{equation}
Let us denote the following event for a given $u\in U$: 
$$A_{6,u}~:~|N^{U}_\leq(u)| > 2\cdot 10^{-3}d_U(u).$$
For every fixed $i\in \{0,1,\ldots,10^3-1\}$ we shall first bound the probability of $A_{6,u}$ under the condition that $i_u=i$. Consider any $u'\in N_U(u)$ such that $0.5d(u)\leq d(u')\leq d(u)$. 
Then, by~(\ref{JuSize}) and~(\ref{3TimesJuLength}), the distance between the right end of any interval that might be assigned as $J(u')$ to $u'$ (say for $i_{u'}=j$) and the left end of the consecutive such interval (for $i_{u'}=j+1$) is larger than:
$$9.9\cdot 10^{-4}d(u') - 3.23 \cdot 10^{-4}d(u') > 6.6 \cdot 10^{-4}d(u') \geq 3.3 \cdot 10^{-4}d(u) > |J(u)|,$$
and hence we may have $J(u')\cap J(u)\neq \emptyset$ for at most one value of $i_{u'}$. Therefore:
$$\mathbf{Pr}\left(J(u')\cap J(u)\neq \emptyset~|~i_u=i\right) \leq 10^{-3}.$$
Thus, by the Chernoff Bound and~(\ref{Udegree}):
$$\mathbf{Pr}(A_{6,u}~|~i_u=i) < 2e^{-\frac{10^{-3}d_U(u)}{3}} < 2e^{-10^{-8}d(u)} <\frac{1}{e\cdot 3\Delta^2},$$
and hence, by the law of total probability:
\begin{equation}\label{PrA6u}
\mathbf{Pr}(A_{6,u}) < \frac{1}{e\cdot 3\Delta^2}.
\end{equation}
Note that each of the events $A_{4,u}$, $A_{5,w}$, $A_{6,u}$ is mutually independent of all other events of such types except at most $3\Delta^2$ of them. Thus, by~(\ref{PrA4u}), (\ref{PrA5w}), (\ref{PrA6u}) and the Lov\'asz Local Lemma,
with positive probability none of these events appears, and hence $F_U\subseteq F'$ and $i_u$ for $u\in U$ can be chosen so that ~(\ref{FUdegreeInU}), (\ref{FUdegreeInW})  and~(\ref{NULeqInU}) hold for all $u\in U$ and $w\in W$. 
 
 Fix any such $F_U$ and $i_u$, $u\in U$, and denote:
 $$G'=(W,E'):=G[W],$$
 i.e. $E'=E(G[W])$.  Let us assign initial weights to all edges outside $W$ by setting:
 \begin{equation}\label{InitialC1Weights}
 \omega_1(e):=\left\{\begin{array}{lcl}
 2&{\rm if} & e\subset U~~{\rm or}~~e\in F_U;\\
 1&{\rm if} & e\in F\smallsetminus F_U.
 \end{array}\right.
 \end{equation}
 A random mode of choosing weights for the remaining edges shall be discussed in Subsection~\ref{Subsection_Rules_W}. At the same time we may still change the weights of some of the edges in $F_W$ and in $E(U)$ later on.

 \subsection{Random Weights of the Edges Inside $W$ -- Preliminaries}\label{Subsection_Wpreliminaries}

 We next proceed towards choosing weights for all the edges of $G'$, thus for most of the edges of $G$ -- these weights shall not be changed in the later part of the construction. 

For this aim we associate with every vertex $v\in W$ an independent random variable $X_v$ valued in $[1.1,2.9]$ with density of the probability distribution:
\begin{equation}\label{DensityGDef}
g(x):=\frac{1}{\ln2.9-\ln1.1}\frac{1}{x}=\frac{1}{\ln\frac{2.9}{1.1}}\frac{1}{x}.
\end{equation}
Further, we associate with every edge $e\in E'$ an independent random variable $X_e\sim U(0,1)$ with uniform distribution over the interval $[0,1]$.

Set: 
\begin{equation}\label{a1a2Def}
a_1:=\frac{2.9}{\left(\frac{2.9}{1.1}\right)^{0.95}},~~~~~~~~~~~~
a_2:=\frac{2.9}{\left(\frac{2.9}{1.1}\right)^{0.45}}
\end{equation} 
(where $0<a_1<a_2<1.9$) and let us define a function $r:[1.1,1.9]\to\mathbb{R}$ where 
$$r(x):=\left\{\begin{array}{lcl}
\frac{x-1}{2}&~{\rm for}~&x\in[1.1,a_1);\\
\frac{x-1}{2}-\frac{1}{\ln\frac{2.9}{1.1}}\cdot\ln\frac{2.9}{1+2\cdot\frac{\ln\frac{2.9}{x}}{\ln\frac{2.9}{1.1}}}&~{\rm for}~&x\in[a_1,a_2];\\
\frac{x-1}{2}-\frac{\ln\frac{2.9}{1.9}}{\ln\frac{2.9}{1.1}}&~{\rm for}~&x\in(a_2,1.9].
\end{array}\right.$$
This function shall be essential while laying down rules of a random choice of weights $1$ and $3$ for edges $uv\in E(W)$ with $X_u$ and $X_v$ relatively small, cf.~(\ref{InitialC1Weights2}). We note that in fact for every $x\in [1.1,1.9]$:
\begin{equation}\label{rValuesBounds}
0<r(x)<0.08.
\end{equation}
This is straightforward to verify for $x\in [1.1,a_1)\cup (a_2,1.9]$. Let us thus consider $x\in [a_1,a_2]$. 
Then~(\ref{rValuesBounds}) follows from the facts that $r(a_1)<0.8$, $r(a_2)>0$, and 
$$
r'(x) = \frac{1}{2} +\frac{1}{\ln \frac{2.9}{1.1}} \cdot \frac{1}{1+2\cdot \frac{\ln\frac{2.9}{x}}{\ln\frac{2.9}{1.1}}} 
\cdot \frac{2}{\ln\frac{2.9}{1.1}} \cdot \frac{-1}{x} < 0
$$ 
inside $[a_1,a_2]$, as this is equivalent (after substituting $z=\frac{2.9}{x}$, with $\frac{2.9}{a_2} = \left(\frac{2.9}{1.1}\right)^{0.45}$ and $\frac{2.9}{a_1} = \left(\frac{2.9}{1.1}\right)^{0.95}$) to the fact that 
$$h(z): = z - \frac{2.9}{4}\cdot \ln\frac{2.9}{1.1}\cdot \left(\ln\frac{2.9}{1.1} +2 \ln z\right) > 0$$
for $z\in \left[\left(\frac{2.9}{1.1}\right)^{0.45}, \left(\frac{2.9}{1.1}\right)^{0.95}\right]$, what in turn holds since 
$h(\left(\frac{2.9}{1.1}\right)^{0.45}) > 0$ and $h'(z) = 1 - \frac{2.9}{2}\cdot \ln\frac{2.9}{1.1} \cdot \frac{1}{z} > 0$ inside $\left[\left(\frac{2.9}{1.1}\right)^{0.45}, \left(\frac{2.9}{1.1}\right)^{0.95}\right]$.

We also define a constant
\begin{equation}\label{OverlineDDef}
\overline{d}:=\int_{1.1}^{1.9}r(x)g(x)dx.
\end{equation}
In order to set down $\overline{d}$ we shall use a substitution:
$$y=1+2\cdot\frac{\ln\frac{2.9}{x}}{\ln\frac{2.9}{1.1}} = 1+\frac{2}{\ln\frac{2.9}{1.1}}\left(\ln2.9-\ln x\right),$$
where:
$$dy=-\frac{2}{x\ln\frac{2.9}{1.1}}dx,$$
$$1+2\cdot\frac{\ln\frac{2.9}{a_1}}{\ln\frac{2.9}{1.1}}  = 1+2\cdot\log_{\frac{2.9}{1.1}}\frac{2.9}{a_1} = 2.9,$$
$$1+2\cdot\frac{\ln\frac{2.9}{a_2}}{\ln\frac{2.9}{1.1}}  = 1+2\cdot\log_{\frac{2.9}{1.1}}\frac{2.9}{a_2} = 1.9.$$
We also note that:
$$\int\ln\frac{2.9}{y} dy = y\ln\frac{2.9}{y}-\int\left(-\frac{1}{y}\right) y ~dy =   y\ln\frac{2.9}{y} + y + const.$$
Thus:
\begin{eqnarray}
\overline{d} &=& \int_{1.1}^{a_1}g(x)\frac{x-1}{2}dx 
+\int_{a_1}^{a_2}g(x)\left(\frac{x-1}{2}-\frac{1}{\ln\frac{2.9}{1.1}}\cdot\ln\frac{2.9}{1+2\cdot\frac{\ln\frac{2.9}{x}}{\ln\frac{2.9}{1.1}}}\right)dx
+\int_{a_2}^{1.9}g(x)\left(\frac{x-1}{2}-\frac{\ln\frac{2.9}{1.9}}{\ln\frac{2.9}{1.1}}\right)dx\nonumber\\
&=&        \int_{1.1}^{1.9}g(x)\frac{x-1}{2}dx 
- \frac{1}{\ln\frac{2.9}{1.1}} \int_{a_1}^{a_2}g(x)\ln\frac{2.9}{1+2\cdot\frac{\ln\frac{2.9}{x}}{\ln\frac{2.9}{1.1}}}dx
-\frac{\ln\frac{2.9}{1.9}}{\ln\frac{2.9}{1.1}}\int_{a_2}^{1.9}g(x)dx \nonumber\\
&=& \frac{1}{2\ln\frac{2.9}{1.1}}\int_{1.1}^{1.9}\left(1-\frac{1}{x}\right)dx 
+ \frac{1}{2\ln\frac{2.9}{1.1}} \int_{a_1}^{a_2}\ln\frac{2.9}{1+2\cdot\frac{\ln\frac{2.9}{x}}{\ln\frac{2.9}{1.1}}}\left(-\frac{2}{\ln\frac{2.9}{1.1}} \frac{1}{x} \right)dx
-\frac{\ln\frac{2.9}{1.9}}{\left(\ln\frac{2.9}{1.1}\right)^2}\int_{a_2}^{1.9}\frac{1}{x}dx \nonumber\\
&=& \frac{1}{2\ln\frac{2.9}{1.1}}    \left(1.9-1.1-\ln\frac{1.9}{1.1}\right) 
- \frac{1}{2\ln\frac{2.9}{1.1}} \int_{1.9}^{2.9}\ln\frac{2.9}{y} dy
-\frac{\ln\frac{2.9}{1.9}}{\left(\ln\frac{2.9}{1.1}\right)^2} \ln\frac{1.9}{a_2} \nonumber\\
&=& \frac{0.8-\ln\frac{1.9}{1.1}}{2\ln\frac{2.9}{1.1}}     
- \frac{1}{2\ln\frac{2.9}{1.1}} \left( 2.9\ln\frac{2.9}{2.9}+2.9- 1.9\ln\frac{2.9}{1.9}-1.9\right)
-\frac{\ln\frac{2.9}{1.9}}{\left(\ln\frac{2.9}{1.1}\right)^2} \left(0.45 \ln\frac{2.9}{1.1} - \ln\frac{2.9}{1.9} \right) \nonumber\\
&=& \frac{1}{\ln\frac{2.9}{1.1}}\left[ -0.1-0.5\left(\ln\frac{1.9}{1.1}+\ln\frac{2.9}{1.9}\right) +\ln\frac{2.9}{1.9} + \frac{\left(\ln\frac{2.9}{1.9}\right)^2}{\ln\frac{2.9}{1.1}}\right] \nonumber\\
&=& \left(\frac{\ln\frac{2.9}{1.9}}{\ln\frac{2.9}{1.1}}\right)^2+\frac{\ln\frac{2.9}{1.9}}{\ln\frac{2.9}{1.1}}+0.25-0.75-\frac{0.1}{\ln\frac{2.9}{1.1}} \nonumber\\
&=& \left(\frac{\ln\frac{2.9}{1.9}}{\ln\frac{2.9}{1.1}}+0.5\right)^2 -\frac{0.1}{\ln\frac{2.9}{1.1}}-0.75\\
&\approx& 0.023. \label{dBarApproximation}
\end{eqnarray}
(More precisely, $\overline{d}\in(0.023,0.024)$.)
Therefore, by~(\ref{rValuesBounds}) and~(\ref{dBarApproximation}), for any $x,y\in[1.1,1.9]$:
\begin{equation}\label{r2-d-relation}
0<\frac{r(x)r(y)}{\overline{d}}< 1.
\end{equation}

\subsection{Weighting Rules within $W$}\label{Subsection_Rules_W}

We may now finally lay down rules defining a $\{1,3\}$-weighting of the edges of $G'$ resulting from our random process. 

We complete the edge weighting $\omega_1$ of $G$, defined in~(\ref{InitialC1Weights}) for the edges outside $W$, by setting for every edge $uv\in E'$ with $X_u\leq X_v$ (cf. in particular~(\ref{r2-d-relation})):
\begin{equation}\label{InitialC1Weights2}
\omega_1(uv):=\left\{\begin{array}{lcl}
3&~~{\rm if}~~& X_v\geq 1.9~~{\rm and}~~X_u\geq \frac{2.9}{(\frac{2.9}{1.1})^\frac{X_v-1}{2}}\\
&& {\rm or}~~X_v<1.9~~{\rm and}~~X_{uv}\leq \frac{r(X_u)\cdot r(X_v)}{\overline{d}};\\
1 &~~{\rm if}~~& X_v\geq 1.9~~{\rm and}~~X_u< \frac{2.9}{(\frac{2.9}{1.1})^\frac{X_v-1}{2}}\\
&& {\rm or}~~X_v<1.9~~{\rm and}~~X_{uv}> \frac{r(X_u)\cdot r(X_v)}{\overline{d}}.
\end{array}
\right.
\end{equation}

Note in particular that by the first condition above, every edge $uv\in E'$ with $X_u\geq 1.9$ and $X_v\geq 1.9$ is weighted $3$, as $\frac{2.9}{(\frac{2.9}{1.1})^\frac{X_v-1}{2}}\leq 1.9$ whenever $X_v\geq 1.9$. 

Thus,  for any edge $uv$, if $X_v\geq 1.9$, then $uv$ is weighted $3$ if and only if $X_u\geq \frac{2.9}{(\frac{2.9}{1.1})^\frac{X_v-1}{2}}$. 

Let $N'_i(v)$ denote the set of edges incident with $v$ in $G'$ and weighted $i$, and set 
$$d'_i(v):=|N'_i(v)|$$ 
for $i=1,3$. 

The resulting weighting $\omega_1$ defines \emph{initial sums} of all vertices $v\in V$, denoted by $s_1(v)$. By~(\ref{InitialC1Weights}) and~(\ref{InitialC1Weights2}) we thus have for every $v\in W$:
\begin{eqnarray}
s_1(v) &=& 2\cdot d_{F_U}(v)+1\cdot \left(d_F(v)-d_{F_U}(v)\right)+3\cdot d'_3(v)+1\cdot d'_1(v) \nonumber\\
&=&d_U(v)+d_{F_U}(v) +d_W(v)+2d'_3(v). \label{S1v}
\end{eqnarray}

\subsection{Near Locations of Initial Sums in $W$}

We shall require two major features from the constructed weighting of the edges. First of all we shall guarantee that
the initial sum of every $v\in W$ is roughly close to $X_v\cdot d(v)$,  more precisely that: 
\begin{equation}\label{s1NearLocation}
s_1(v)\in\left[d_U(v)+d_{F_U}(v) +X_vd_W(v)-10^{-9}d_W(v),d_U(v)+d_{F_U}(v) +X_vd_W(v)+10^{-9}d_W(v)\right].
\end{equation}
We thus define the following event for every $v\in W$:
$$
A_v:~ s_1(v)\not\in\left[d_U(v)+d_{F_U}(v) +X_vd_W(v)-10^{-9}d_W(v),d_U(v)+d_{F_U}(v) +X_vd_W(v)+10^{-9}d_W(v)\right].
$$
By~(\ref{S1v}), this is equivalent to the following event: 
$$A_v:~ d'_3(v)\not\in\left[\frac{X_v-1}{2}d_W(v)-\frac{10^{-9}}{2}d_W(v),\frac{X_v-1}{2}d_W(v)+\frac{10^{-9}}{2}d_W(v)
\right].$$
For any given $\alpha\in[1.1,2.9]$, we shall first bound the conditional probability:
$$\mathbf{Pr}\left(A_v~|~X_v=\alpha\right).$$
For each fixed $\alpha\in[1.1,2.9]$ this probability equals to the probability that $A_v$ holds but in a slightly altered experiment where the (new) random variable $X_v$ equals $\alpha$ with pro\-ba\-bi\-li\-ty $1$ (and the rest of variables remain unchanged). For distinction, the probability of a given event with respect to this new experiment shall be denoted by $\mathbf{Pr}_\alpha$, hence:
$$\mathbf{Pr}_\alpha\left(A_v\right)=\mathbf{Pr}\left(A_v~|~X_v=\alpha\right).$$

For any $u\in N_W(v)$, denote the event:
$$A_{v,u}: ~uv\in N'_3(v).$$
(Note that due to $\alpha$ being fixed as the value of $X_v$, the events $A_{v,u}$ are independent for distinct $u\in N_W(v)$.)\\

Assume first that $\alpha\geq 1.9$. (Note that, as $\alpha \leq 2.9$, we have $\frac{2.9}{(\frac{2.9}{1.1})^\frac{\alpha-1}{2}} > 1.1$ then.)

Then for any fixed $u\in N_W(v)$, by~(\ref{InitialC1Weights2}) and~(\ref{DensityGDef}): 
\begin{eqnarray}
\mathbf{Pr}_\alpha\left(A_{v,u} \right) &=& \mathbf{Pr}\left(X_u\geq \frac{2.9}{(\frac{2.9}{1.1})^\frac{\alpha-1}{2}}\right)
= \int\limits_{\frac{2.9}{(\frac{2.9}{1.1})^\frac{\alpha-1}{2}}}^{2.9}\frac{1}{\ln\frac{2.9}{1.1}}\frac{1}{x}dx\nonumber\\
&=& \frac{1}{\ln\frac{2.9}{1.1}} \left(\ln2.9 - \ln\left(\frac{2.9}{(\frac{2.9}{1.1})^\frac{\alpha-1}{2}}\right)\right)
= \frac{1}{\ln\frac{2.9}{1.1}}\cdot \ln \left(\left(\frac{2.9}{1.1}\right)^\frac{\alpha-1}{2}\right) \nonumber\\
&=& \frac{\alpha-1}{2}. \label{Buv1case}
\end{eqnarray}

Assume now that $\alpha < 1.9$. Then for any fixed $u\in N_W(v)$: 
\begin{eqnarray}
\mathbf{Pr}_\alpha\left(A_{v,u} \right) &=& \mathbf{Pr}_\alpha\left(uv\in N'_3(v) \wedge X_u\geq 1.9\right) + \mathbf{Pr}_\alpha\left(uv\in N'_3(v) \wedge X_u< 1.9\right)\nonumber\\
&=& \mathbf{Pr}_\alpha\left(uv\in N'_3(v) \wedge X_u\geq 1.9\right) + \mathbf{Pr}\left(X_{uv}\leq  \frac{r(X_u)\cdot r(\alpha)}{\overline{d}} \wedge X_u< 1.9\right).\nonumber\\ \label{Buv2caseA}
\end{eqnarray}
Moreover, for any fixed $\beta\geq 1.9$,
\begin{eqnarray}
\mathbf{Pr}\left(X_{uv}\leq  \frac{r(X_u)\cdot r(\alpha)}{\overline{d}} \wedge X_u< 1.9~|~X_u=\beta \right) &=& 0, \label{Buv2caseB}
\end{eqnarray}
while for a fixed $\beta< 1.9$,
\begin{eqnarray}
\mathbf{Pr}\left(X_{uv}\leq  \frac{r(X_u)\cdot r(\alpha)}{\overline{d}} \wedge X_u< 1.9~|~X_u=\beta \right) &=& \mathbf{Pr}\left(X_{uv}\leq  \frac{r(\beta)\cdot r(\alpha)}{\overline{d}} \right) \nonumber\\
&=& \frac{r(\beta)\cdot r(\alpha)}{\overline{d}}.
\label{Buv2caseC}
\end{eqnarray}
Hence, by~(\ref{Buv2caseA}), (\ref{Buv2caseB}), (\ref{Buv2caseC}) and~(\ref{OverlineDDef}),
\begin{eqnarray}
\mathbf{Pr}_\alpha\left(A_{v,u} \right) &=& \mathbf{Pr}_\alpha\left(uv\in N'_3(v) \wedge X_u\geq 1.9\right) +  \int\limits_{1.1}^{1.9} \frac{r(x)\cdot r(\alpha)}{\overline{d}}g(x)dx \nonumber\\
&=& \mathbf{Pr}_\alpha\left(uv\in N'_3(v) \wedge X_u\geq 1.9\right) +  \frac{r(\alpha)}{\overline{d}} \int\limits_{1.1}^{1.9} r(x)g(x)dx \nonumber\\
&=& \mathbf{Pr}_\alpha\left(uv\in N'_3(v) \wedge X_u\geq 1.9\right) + r(\alpha). \label{Buv2caseD}
\end{eqnarray}
Therefore, by~(\ref{Buv2caseD}), (\ref{InitialC1Weights2}) and~(\ref{a1a2Def}), for $a_2 < \alpha< 1.9$ (as $\frac{2.9}{(\frac{2.9}{1.1})^\frac{X_u-1}{2}}\leq a_2$ for $X_u\geq 1.9$),
\begin{eqnarray}
\mathbf{Pr}_\alpha\left(A_{v,u} \right) &=&  \mathbf{Pr}\left( X_u\geq 1.9\right)+ r(\alpha)
= \int\limits_{1.9}^{2.9}\frac{1}{\ln\frac{2.9}{1.1}}\frac{1}{x}dx + \left(\frac{\alpha-1}{2}-\frac{\ln\frac{2.9}{1.9}}{\ln\frac{2.9}{1.1}}\right)  \nonumber\\
&=& \frac{\ln\frac{2.9}{1.9}}{\ln\frac{2.9}{1.1}}  + \frac{\alpha-1}{2}-\frac{\ln\frac{2.9}{1.9}}{\ln\frac{2.9}{1.1}} 
= \frac{\alpha-1}{2}, \label{Buv2case1subcase}
\end{eqnarray}
while for $a_1\leq \alpha \leq a_2$ (as then: $1.9\leq 1+2\frac{\ln\frac{2.9}{\alpha}}{\ln\frac{2.9}{1.1}} \leq 2.9$), 
\begin{eqnarray}
\mathbf{Pr}_\alpha\left(A_{v,u} \right) 
 &=&  \mathbf{Pr}\left( \alpha \geq \frac{2.9}{(\frac{2.9}{1.1})^\frac{X_u-1}{2}} \wedge X_u\geq 1.9\right)+ r(\alpha)\nonumber\\
  &=&  \mathbf{Pr}\left( X_u \geq 1+2\frac{\ln\frac{2.9}{\alpha}}{\ln\frac{2.9}{1.1}} \wedge X_u\geq 1.9\right)+ r(\alpha) \nonumber\\
    &=&  \mathbf{Pr}\left( X_u \geq 1+2\frac{\ln\frac{2.9}{\alpha}}{\ln\frac{2.9}{1.1}} \right)+ r(\alpha) = \int\limits_{1+2\frac{\ln\frac{2.9}{\alpha}}{\ln\frac{2.9}{1.1}}}^{2.9}\frac{1}{\ln\frac{2.9}{1.1}}\frac{1}{x}dx + r(\alpha)  \nonumber\\
&=& \left(\frac{1}{\ln\frac{2.9}{1.1}}\cdot\ln\frac{2.9}{1+2\cdot\frac{\ln\frac{2.9}{\alpha}}{\ln\frac{2.9}{1.1}}}\right) + \left(\frac{\alpha-1}{2}-\frac{1}{\ln\frac{2.9}{1.1}}\cdot\ln\frac{2.9}{1+2\cdot\frac{\ln\frac{2.9}{\alpha}}{\ln\frac{2.9}{1.1}}}\right) \nonumber\\ 
&=& \frac{\alpha-1}{2}, \label{Buv2case2subcase}
\end{eqnarray}
and finally, for $1.1\leq \alpha < a_1$ (as 
$a_1\leq \frac{2.9}{(\frac{2.9}{1.1})^\frac{X_u-1}{2}}$), 
\begin{eqnarray}
\mathbf{Pr}_\alpha\left(A_{v,u} \right) 
 &=& 0+ r(\alpha) ~=~ \frac{\alpha-1}{2}. \label{Buv2case3subcase}
\end{eqnarray}
Since all the  events $A_{v,u}$ are independent for distinct $u\in N_W(v)$ (in our modified probability space, where $X_v=\alpha$), by (\ref{Buv1case}), (\ref{Buv2case1subcase}), (\ref{Buv2case2subcase}), (\ref{Buv2case3subcase}), the Chernoff Bound and~(\ref{Udegree}),
$$\mathbf{Pr}_\alpha(A_v) 
< 2e^{-\frac{\left(\frac{10^{-9}}{2}d_W(v)\right)^2}{3\frac{\alpha-1}{2}d_W(v)}}
= 2e^{-\frac{10^{-18}d_W(v)}{6(\alpha-1)}}
\leq 2e^{-\frac{10^{-18}0.9d(v)}{11.4}}
\leq 2e^{-\frac{90\ln\Delta}{11.4}}
\leq\frac{1}{2e\Delta^2}.$$  
Hence, in the original probability space:
$$\mathbf{Pr}\left(A_v~|~X_v=\alpha\right) < \frac{1}{2e\Delta^2},$$ 
and thus:
\begin{equation}\label{FinalPrAvUpeerBound}
\mathbf{Pr}\left(A_v\right) 
= \int_{1.1}^{2.9}g(\alpha)\mathbf{Pr}\left(A_v~|~X_v=\alpha\right)d\alpha 
\leq \int_{1.1}^{2.9}g(\alpha) \frac{1}{2e\Delta^2}d\alpha =  \frac{1}{2e\Delta^2} . 
\end{equation}

\subsection{Distribution of Near Locations of Final Sums in $W$}

The final sums of the vertices in $W$, resulting from increasing weights of some of the edges in $F_W$, shall be a bit bigger than the initial ones. Thereby we shall be able to control their precise values. 
The second required feature (cf.~(\ref{NWsamllerCondition})) of our major random process shall thus concern well distribution of so-called near-locations of sums in $W$ -- auxiliary quantities, which are slightly larger than initial sums, but shall eventually be very close to the final sums in $W$.

For every $v\in W$, we define: 
\begin{eqnarray}
l(v) &:=& 2^{\left\lfloor\log_2\left(10^{-9}d_W(v)\right)\right\rfloor}, \nonumber\\
\mathcal{I}(v) &:=& \left\{\left[k\cdot l(v),(k+1)\cdot l(v)\right):k\in\mathbb{Z}\right\}, \nonumber
\end{eqnarray}
and note that the length $l(v)$ of every interval in $\mathcal{I}(v)$ fulfills:
\begin{equation}\label{lvEstimations}
0.5\cdot10^{-9}d_W(v)\leq l(v)\leq 10^{-9}d_W(v).
\end{equation}
We also define the following quantity, which we shall call the \emph{near location of the final sum} of $v$:
\begin{equation}\label{s0Def}
s_0(v):=d_U(v)+d_{F_U}(v) +X_vd_W(v)+3l(v).
\end{equation}
We shall later guarantee that the final sum of $v$ shall belong to exactly the same (short) interval from the family 
$\mathcal{I}(v)$ as $s_0(v)$. Let thus $I(v)=[i_0(v),i_1(v))\in \mathcal{I}(v)$ be such an interval (of length $l(v)$) that 
\begin{equation}\label{S0vBelongingReq}
s_0(v)\in I(v).
\end{equation}
Denote:
$$i'_0(v)=\frac{i_0(v)-d_U(v)-d_{F_U}(v)-3l(v)}{d_W(v)},$$
$$i'_1(v)=\frac{i_1(v)-d_U(v)-d_{F_U}(v)-3l(v)}{d_W(v)}.$$
Then:
$$X_v\in \left[i'_0(v),i'_1(v)\right),$$
and thus, by~(\ref{lvEstimations}):
\begin{equation}\label{i'0vLowerBound}
i'_0(v)= i'_1(v)-\frac{l(v)}{d_W(v)} > 1.1 - 10^{-9}.
\end{equation}
Let 
\begin{equation}\label{NWLeqDefinition}
N^W_\leq(v):=\left\{u\in N_W(v): d_W(u)\leq d_W(v)\right\}.
\end{equation}
We shall guarantee that for every $v\in W$:
\begin{equation}\label{NWsamllerCondition}
\left|\left\{u\in N^W_\leq(v) ~|~ s_0(u)\in I(v) \right\}\right| \leq 0.95 \cdot l(v).
\end{equation}
Let us thus consider the following event for any given $v\in W$:
$$A'_v: \left|\left\{u\in N^W_\leq(v) ~|~ s_0(u)\in I(v) \right\}\right| > 0.95 \cdot l(v).$$

For any fixed $\alpha\in[1.1,2.9]$, we shall first bound the conditional probability:
$$\mathbf{Pr}\left(A'_v~|~X_v=\alpha\right).$$
Analogously as above, for the given $\alpha\in[1.1,2.9]$ this is equal to the probability that $A'_v$ holds but in a new experiment where the (new) random variable $X_v$ equals $\alpha$ with probability $1$. We recall that for distinction, the probability of a given event with respect to this new experiment shall be denoted by $\mathbf{Pr}_\alpha$, hence: $\mathbf{Pr}_\alpha\left(A'_v\right)=\mathbf{Pr}\left(A'_v~|~X_v=\alpha\right)$.

We start from bounding the probability (within the modified model, with $X_v=\alpha$) of the following event for a fixed $u\in N^W_\leq(v)$:
$$A'_{v,u}:~ s_0(u)\in I(v).$$
(Note the events $A'_{v,u}$ are independent for distinct $u\in N^W_\leq (v)$ then.)
 
For this aim however we need to show beforehand  that:  
\begin{equation}\label{dUdFU-difference}
d_U(v)-d_U(u)+d_{F_U}(v)-d_{F_U}(u) \geq -5\cdot 10^{-6} d_W(v).
\end{equation}
First note that as $d_W(u)\leq d_W(v)$, by~(\ref{Udegree}) we obtain that $d(u)(1-10^{-4}-10^{-6}) \leq d(v)(1-10^{-4}+10^{-6})$, and hence:
\begin{equation}\label{dudv-relation_in_W}
d(u) \leq \frac{1-10^{-4}+10^{-6}}{1-10^{-4}-10^{-6}} d(v) < (1+3\cdot10^{-6})d(v).
\end{equation}
Further, by~(\ref{FUdegreeInW}), (\ref{F'_definition}), (\ref{F_WdegreeInW}), (\ref{Udegree}) and~(\ref{dudv-relation_in_W}),
\begin{eqnarray}
d_U(v)+d_{F_U}(v) &\geq& d_U(v) + \left(\frac{1-10^{-3}}{2}-10^{-5}\right)d_{F'}(v) \nonumber\\
&\geq& d_U(v) + \left(\frac{1-10^{-3}}{2}-10^{-5}\right)(1-10^{-4}-10^{-6}) d_U(v)  \nonumber \\
&\geq& \left[1 + \left(\frac{1-10^{-3}}{2}-10^{-5}\right)(1-10^{-4}-10^{-6}) \right](10^{-4}-10^{-6}) d(v)\nonumber\\
&\geq& 0.0001484 d(v),\label{dUvdFUvsumInW}
\end{eqnarray}
\begin{eqnarray}
d_U(u)+d_{F_U}(u) &\leq& d_U(v) + \left(\frac{1-10^{-3}}{2}+10^{-5}\right)d_{F'}(v) \nonumber\\
&\leq& d_U(v) + \left(\frac{1-10^{-3}}{2}+10^{-5}\right)(1-10^{-4}+10^{-6}) d_U(v)  \nonumber\\
&\leq& \left[1 + \left(\frac{1-10^{-3}}{2}+10^{-5}\right)(1-10^{-4}+10^{-6})\right] (10^{-4}+10^{-6})d(u) \nonumber\\
&\leq& \left[1 + \left(\frac{1-10^{-3}}{2}+10^{-5}\right)(1-10^{-4}+10^{-6})\right] (10^{-4}+10^{-6})(1+3\cdot10^{-6})d(v) \nonumber\\
&\leq& 0.0001515 d(v)\label{dUudFUusumInW}.
\end{eqnarray}
Thus, by~({\ref{dUvdFUvsumInW}) and~(\ref{dUudFUusumInW}), we have: $d_U(v)-d_U(u)+d_{F_U}(v)-d_{F_U}(u) \geq -4\cdot 10^{-6} d(v)$, and hence (\ref{dUdFU-difference}) follows by~(\ref{Udegree}).

Denote temporarily: 
$$i''_0=\frac{i_0(v)-d_U(u)-d_{F_U}(u)-3l(u)}{d_W(u)},$$
$$i''_1=\frac{i_1(v)-d_U(u)-d_{F_U}(u)-3l(u)}{d_W(u)}.$$
Then, as $\ln(1+x)\leq x$, by the definitions of $s_0(u)$ and $I(v)$, (\ref{i'0vLowerBound}) and~(\ref{dUdFU-difference}):
\begin{eqnarray}
\mathbf{Pr}_\alpha\left(A'_{v,u}\right) &=& \mathbf{Pr}_\alpha\left(X_u \in [i''_0,i''_1)\right)  
=  \mathbf{Pr}_\alpha\left(X_u \in [i''_0,i''_1) \cap [1.1,2.9]\right)  \nonumber \\
&=& \frac{1}{\ln\frac{2.9}{1.1}}\int_{[i''_0,i''_1] \cap [1.1,2.9]} \frac{1}{x}dx
\leq \frac{1}{\ln\frac{2.9}{1.1}}\int_{[i''_0,i''_1]} \frac{1}{x}dx \nonumber\\
&=&  \frac{1}{\ln\frac{2.9}{1.1}}  \ln\frac{i''_1}{i''_0} 
= \frac{1}{\ln\frac{2.9}{1.1}}  \ln\frac{d_W(u)i''_1}{d_W(u)i''_0} 
=  \frac{1}{\ln\frac{2.9}{1.1}}  \ln\frac{d_W(u)i''_0+l(v)}{d_W(u)i''_0} \nonumber\\
&\leq& \frac{1}{\ln\frac{2.9}{1.1}} \frac{l(v)}{d_W(u)i''_0} 
=  \frac{1}{\ln\frac{2.9}{1.1}} \frac{l(v)}{d_W(v)} \cdot \frac{1}{\frac{i_0(v)-d_U(u)-d_{F_U}(u)-3l(u)}{d_W(v)}} \nonumber\\
&=& \frac{1}{\ln\frac{2.9}{1.1}} \frac{l(v)}{d_W(v)}  \cdot \frac{1}{\frac{i'_0(v)d_W(v)+d_U(v)-d_U(u)+d_{F_U}(v)-d_{F_U}(u)+3l(v)-3l(u)}{d_W(v)}} \nonumber\\
&\leq& \frac{1}{\ln\frac{2.9}{1.1}} \frac{l(v)}{d_W(v)} \cdot \frac{1}{1.1 - 10^{-9}+\frac{d_U(v)-d_U(u)+d_{F_U}(v)-d_{F_U}(u)}{d_W(v)}} \nonumber\\
&\leq&  \frac{l(v)}{d_W(v)} \cdot \frac{1}{\ln\frac{2.9}{1.1}} \cdot \frac{1}{1.1 - 6\cdot 10^{-6}} 
\leq 0.94  \frac{l(v)}{d_W(v)}. \label{PrA'vuUperBound}
\end{eqnarray}

As the  vents $A'_{v,u}$ are independent for distinct $u\in N^W_\leq(v)$ (in our $\alpha$ modified probability space),
by (\ref{PrA'vuUperBound}), (\ref{lvEstimations}), (\ref{Udegree}), the Chernoff Bound and the fact that $|N^W_\leq(v)|\leq d_W(v)$,
$$\mathbf{Pr}_\alpha(A'_v) < 2e^{-\frac{\left(0.01\cdot l(v)\right)^2}{3 \cdot0.94\cdot l(v)}} 
=  2e^{-\frac{10^{-4}\cdot l(v)}{2.82}} 
\leq 2e^{-\frac{10^{-13}\cdot 0.5d_W(v)}{2.82}} 
\leq 2e^{-10^{-14}d(v)}  
\leq 2e^{-10^{6}\ln\Delta}  
\leq \frac{1}{2e\Delta^2}.$$ 
Hence, in the original probability space:
$$\mathbf{Pr}\left(A'_v~|~X_v=\alpha\right) < \frac{1}{2e\Delta^2},$$ 
and thus: 
\begin{equation}\label{FinalPrA'vUpeerBound}
\mathbf{Pr}\left(A'_v\right) = \int_{1.1}^{2.9}g(\alpha)\mathbf{Pr}\left(A'_v~|~X_v=\alpha\right)d\alpha 
\leq \int_{1.1}^{2.9}g(\alpha) \frac{1}{2e\Delta^2}d\alpha = \frac{1}{2e\Delta^2}. 
\end{equation}

 \subsection{Setting Final Sums in $W$}\label{Subsection_Wfinal}

As each event $A_v$ and each event $A'_v$ is mutually independent of all other events of such types except possibly those associated with vertices at distance at most $2$ from $v$ in $G'$, hence less than $2\Delta^2$ such events, by (\ref{FinalPrAvUpeerBound}), (\ref{FinalPrA'vUpeerBound}) and the Lov\'asz Local Lemma, the values of $X_v$  can be chosen so that none of the events $A_v$, $A'_v$, $v\in W$ appears, i.e. such that~(\ref{s1NearLocation}) and~(\ref{NWsamllerCondition}) are fulfilled. Then, by~(\ref{s1NearLocation}), (\ref{s0Def}), (\ref{lvEstimations}) and~(\ref{S0vBelongingReq}), for every $v\in W$:
\begin{eqnarray}
s_1(v) &\leq & d_U(v)+d_{F_U}(v) + (X_v+10^{-9})d_W(v) 
= s_0(v)-3l(v)+10^{-9}d_W(v) \nonumber\\
&\leq& s_0(v)-l(v) \leq i_0(v), \label{s1vi0vLowerBound}
\end{eqnarray}
\begin{eqnarray}
s_1(v) &\geq & d_U(v)+d_{F_U}(v) + (X_v-10^{-9})d_W(v) 
= s_0(v)-3l(v)-10^{-9}d_W(v) \nonumber\\
&\geq& s_0(v) -5l(v) 
\geq i_1(v)-6l(v). \label{s1vi1vUpperBound}
\end{eqnarray}

We further note that by~ (\ref{Udegree}), (\ref{F_WdegreeInW})  and~(\ref{lvEstimations}), for each $v\in W$:
\begin{equation}\label{dFWLarger}
d_{F_W}(v) > 6l(v).
\end{equation}

We now arbitrarily arrange the vertices in $W$ into a sequence $v_1,v_2,v_3,\ldots$ so that for every $i$, $d_W(v_i)\leq d_W(v_{i+1})$. We then  analyze these vertices one after another consistently with the fixed ordering, and associate with a currently analyzed vertex $v$ an appropriate integer \emph{sum addition} $a(v)$ (we shall eventually increase the sum of $v$ by exactly this quantity) so that
\begin{equation}\label{avReq1}
s_1(v)+a(v) \in I(v)
\end{equation}
(note that by~(\ref{s1vi0vLowerBound}) and~(\ref{s1vi1vUpperBound}) this means we must thus have: $a(v)\in [0,6l(v)]$) and so that
\begin{equation}\label{avReq2}
s_1(v)+a(v) \not\equiv 0,1 ({\rm mod}~100),
\end{equation}
and finally so that for every neighbour $u$ of $v$ in $W$ preceding $v$ in the sequence (thus $u\in N^W_{\leq}(v)$, cf.~(\ref{NWLeqDefinition})), hence with already fixed $a(u)$, we have:
\begin{equation}\label{avReq3}
s_1(v)+a(v)\neq s_1(u)+a(u).
\end{equation}
This is feasible, as by~(\ref{NWsamllerCondition}) we may have only limited number of $u\in N^W_{\leq}(v)$ with  $s_1(u)+a(u)\in I(v)$, because then (since  $s_1(u)+a(u)\in I(u)$) we must have $I(u)\cap I(v)\neq \emptyset$, and thus by the definitions of $I(u)$, $I(v)$ (and $l(v)$, $\mathcal{I}(v)$) and the fact that $d_W(u)\leq d_W(v)$ we must consequently have that $I(u)\subseteq I(v)$, and thus $s_0(u)\in I(u)\subseteq I(v)$, hence by~(\ref{NWsamllerCondition}) there can be at most $0.95l(v)$ such $u$. Thus, as at least $0.97l(v)$ integers in $I(v)$ are not congruent to $0$ nor $1$ modulo $100$, we may obviously choose $a(v)$ fulfilling~(\ref{avReq1}), (\ref{avReq2}) and~(\ref{avReq3}), where by~(\ref{s1vi0vLowerBound}), (\ref{s1vi1vUpperBound}) and~(\ref{avReq1}),
\begin{equation}\label{avAtMost}
a(v)\in [0,6l(v)].
\end{equation}

After finishing analysis of all vertices in $W$ and fixing all their sum additions, for every vertex $v\in W$ we choose  arbitrary $a(v)$ edges from $F_W$ which are incident with $v$ (note that by~(\ref{dFWLarger}) and~(\ref{avAtMost}) there are enough of them), and we change their weights from $1$ to $2$ (cf.~(\ref{InitialC1Weights})). We denote the sum induced by the obtained weighting $\omega_2$ of $G$ by $s_2(v)$ for $v\in V$. Note that then for every $v\in W$, we have:
\begin{equation}\label{S2S1Av}
s_2(v)=s_1(v)+a(v)
\end{equation}
and thus, by~(\ref{avReq3}) and the algorithm above, all neighbours in $W$ are already sum distinguished. We shall not be changing their sums further on.

 \subsection{Algorithm Fixing Sums in $U$}\label{Subsection_Ufinal}

Within the following concluding algorithm, based on ideas introduced by Kalkowski~\cite{Kalkowski12} and developed by Kalkowski, Karo\'nski and Pfender~\cite{KalKarPf_123}, we shall change weights of edges inside $U$ exclusively, denoting the obtained final sums of all the vertices $v\in V$ by $s_3(v)$. We thus note that by~(\ref{FUdegreeInU}), (\ref{F'_definition}), (\ref{WandFDefinition}), (\ref{Udegree}) and~(\ref{3TimesJuLength}), the following shall hold at the end of our construction for every vertex $u\in U$: 
\begin{eqnarray}
  d(u) ~\leq~ s_3(u) &\leq& 2(d_{F_U}(u)+d_{F_W}(u))+(d_W(u)-d_{F_U}(u)-d_{F_W}(u)) +3d_U(u) \label{duBoundFors3uInU} \\
    &=& d(u)+d_{F_U}(u)+d_{F_W}(u)+2d_U(u) \nonumber\\
    &\leq& d(u)+10^{-3}i_ud_{F'}(u) + 10^{-5}d(u) +d_{F_W}(u)+2d_U(u) \label{s3uSmallerRightJu}\\
    &=& 2d(u) + (10^{-3}i_u-1)d_{F'}(u) + 10^{-5}d(u) + d_U(u)\nonumber\\
    &\leq& 2d(u) - 10^{-3} d_{F'}(u) + 10^{-5}d(u) + (10^{-4}+10^{-6})d(u) \nonumber\\
    &<& 2d(u), \label{2duBoundFors3uInU}
\end{eqnarray}
and moreover,
\begin{eqnarray}\label{du2duBoundsFors3uInU}
  s_3(u) &\geq& 2d_{F_U}(u)+(d_W(u)-d_{F_U}(u)) +d_U(u) 
  = d(u)+d_{F_U}(u)\nonumber\\
    &\geq& d(u)+10^{-3}i_ud_{F'}(u) - 10^{-5}d(u). \label{s3uLargerLeftJu}
\end{eqnarray}
Therefore, by~(\ref{duBoundFors3uInU}), (\ref{2duBoundFors3uInU}), (\ref{s3uSmallerRightJu})  and~(\ref{s3uLargerLeftJu}), we shall have (cf.~(\ref{JuDefinition})):
\begin{equation}\label{S3BelongsToIntervals}
s_3(u)\in [d(u),2d(u)]~~~~~~{\rm and}~~~~~~ s_3(u)\in J(u).
\end{equation}

We associate with every $u\in U$ a set of edges $E^*(u)\subseteq (E(U)\cap E(u))$ such that: 
\begin{eqnarray}
E^*(u)\cap E^*(v)= \emptyset~~~&{\rm for}&~~~u,v\in U,~u\neq v,\label{EmptyE*Intersections}\\
|E^*(u)|\geq 0.5 d_U(u)-1~~~&{\rm for}&~~~u\in U.\label{E*uLargeEnoughReq}
\end{eqnarray}
For this aim it is e.g. sufficient to add temporarily to the graph $G[U]$ a new auxiliary vertex $u_0$ (if necessary) and join it by a single edge with every vertex of odd degree in $G[U]$. Each component of the resulting graph shall then be eulerian, and thus we may traverse any fixed Euler tour in each such component, starting form $u$ if it belongs to the component, and direct every traversed edge consistently with our direction of movement along the tour. For every $u\in U$ we then include in $E^*(u)$ each edge in  $E(U)\cap E(u)$ outgoing from $u$. Finally, we remove $u_0$ (and forget of all the directions); it is straightforward to verify that~(\ref{EmptyE*Intersections}) and~(\ref{E*uLargeEnoughReq}) shall be fulfilled then.

We next arrange the vertices of $U$ into a sequence $u_1,u_2,u_3,\ldots$ such that $d(u_i)\leq d(u_{i+1})$ for every $i$. Then we analyze one by one all consecutive vertices in the sequence, and the moment we analyze a given vertex $u\in U$, we associate to it a set $S(u)$ chosen from the family:
\begin{equation}\label{DefinitionMathcalS}
\mathcal{S}:=\{\{100i,100i+1\}:~ i\in \mathbb{Z} \}
\end{equation}
so that 
\begin{equation}\label{SuNeqSvCond}
S(u)\neq S(v) ~~~~{\rm for~ every}~~~~ v\in N^U_\leq(u)
\end{equation}
preceding $u$ in the sequence. We moreover guarantee at the same time that the sum of $u$ belongs to $S(u)$ and never leaves this set later on. For this aim we permit to change the weight of every edge $uv\in E^*(u)$ while analyzing $u$, remembering to assure that the resulting sum at $v$ belongs to $S(v)$ if $v$ precedes $u$ in the sequence.  Note that as initially all weights within $U$ were set to $2$ (cf.~(\ref{InitialC1Weights})) and each weight shall be altered at most once due to~(\ref{EmptyE*Intersections}), we may always change the weight of each such edge by $1$ (keeping it in the range $[1,3]$), hence we may  thereby obtain altogether at least $|E^*(u)|$ distinct sums for $u$, which are consecutive integers. By~(\ref{E*uLargeEnoughReq}) the set of these possible sums must contain more than $4\cdot10^{-3}d_U(u)$ sets from $\mathcal{S}$, while by~(\ref{NULeqInU}) less than half of these can already be assigned to vertices in $N^U_\leq(u)$. Thus we may choose among these sums one, say $s^*$ such that $\{s^*,s^*+1\}\in \mathcal{S}$ and $\{s^*,s^*+1\} \neq S(v)$ for every $v\in N^U_\leq(u)$ preceding $u$ in the sequence. We then perform admissible weight changes of edges in $E^*(u)$ so that the sum at $u$ equals $s^*$, and set $S(u)=\{s^*,s^*+1\}$. After analyzing all vertices in $U$ all our requirements are thus fulfilled, in particular~(\ref{SuNeqSvCond}). We denote the obtained final weighting of the edges of $G$ by $\omega_3$, and the resulting sum at every vertex $v$ in $G$ by $s_3(v)$.

Suppose now that $u_iu_j\in E(U)$ with $i<j$, hence $d(u_i)\leq d(u_j)$. 

If $d(u_i) < 0.5d(u_j)$, then by~(\ref{S3BelongsToIntervals}), $s_3(u_i)\leq 2d(u_i)<d(u_j)\leq s_3(u_j)$.

If on the other hand $d(u_i)\in [0.5d(u_j),d(u_j)]$ and $J(u_i)\cap J(u_j)=\emptyset$, then $s_3(u_i) \neq s_3(u_j)$, as by~(\ref{S3BelongsToIntervals}), $s_3(u_i)\in J(u_i)$ and $s_3(u_j)\in J(u_j)$.

If finally  $d(u_i)\in [0.5d(u_j),d(u_j)]$ and $J(u_i)\cap J(u_j)\neq\emptyset$, then by~(\ref{NULeqDef}), $u_i \in N^U_\leq(u_j)$ and thus, by~(\ref{SuNeqSvCond}), $S(u_i)\neq S(u_j)$, hence $S(u_i)\cap S(u_j)=\emptyset$ (cf.~\ref{DefinitionMathcalS}), and consequently, $s_3(u_i) \neq s_3(u_j)$.

All neighbours in $U$ are thus sum distinguished, similarly as all neighbours in $W$ (as we have not changed their sums in this final stage, cf. the end of the previous subsection). 

To conclude, note that by~(\ref{DefinitionMathcalS}), for every $u\in U$, we have:
$$s_3(u)\equiv 0,1~ ({\rm mod}~100),$$
while by~(\ref{avReq2}) and~(\ref{S2S1Av}), for every $v\in W$:
$$s_3(v) = s_2(v)  \not\equiv 0,1 ({\rm mod}~100).$$
All neighbours in $G$ are thus sum distinguished, hence $\omega_3$ is a desired vertex-colouring $3$-weighting of $G$. \qed


\end{document}